\documentclass{amsart}

\usepackage{amsmath,amssymb,xy}

\xyoption{all}
\CompileMatrices

\theoremstyle{plain}

\newtheorem{theorem}{Theorem}[section]

\newtheorem{lemma}[theorem]{Lemma}
\newtheorem{corollary}[theorem]{Corollary}

\newtheorem{mainthm}{Theorem}

\theoremstyle{definition}

\newtheorem{remark}[theorem]{Remark}

\renewcommand{\O}{\mathcal{O}}
\newcommand{\Q}{\mathbf{Q}}
\newcommand{\Ql}{\Q_{\ell}}
\newcommand{\Qp}{\Q_{p}}
\newcommand{\Z}{\mathbf{Z}}
\newcommand{\Zp}{\Z_{p}}

\newcommand{\Fo}{F}
\newcommand{\Fop}{F'}
\newcommand{\Fopp}{F''}
\newcommand{\Fi}{F_{\infty}}
\newcommand{\Fip}{\Fi'}
\newcommand{\Fipp}{\Fi''}
\newcommand{\Fiw}{F_{\infty,w}}
\newcommand{\Fipwp}{F'_{\infty,w'}}
\newcommand{\Fippwpp}{F''_{\infty,w''}}
\newcommand{\Fov}{F_{v}}
\newcommand{\GFo}{G_{\Fo}}
\newcommand{\GFov}{G_{v}}
\newcommand{\GFopvp}{G_{v'}}
\newcommand{\GF}{G_{F_{\infty}}}

\newcommand{\Gv}{G_{v}}

\newcommand{\Iv}{I_{v}}

\newcommand{\Qbar}{\bar{\mathbf{Q}}}
\newcommand{\Qpbar}{\Qbar_{p}}

\newcommand{\Ac}{A_{\chi}}

\newcommand{\crit}{\text{cr}}
\newcommand{\inj}{\hookrightarrow}
\newcommand{\surj}{\twoheadrightarrow}
\renewcommand{\th}{\text{th}}

\newcommand{\Qinf}{\Q_\infty}
\newcommand{\Qinfv}{\Q_{\infty,v}}

\DeclareMathOperator{\alg}{alg}
\DeclareMathOperator{\an}{an}
\DeclareMathOperator{\Gal}{Gal}
\DeclareMathOperator{\GL}{GL}
\DeclareMathOperator{\Hom}{Hom}
\DeclareMathOperator{\im}{im}
\DeclareMathOperator{\Ind}{Ind}

\DeclareMathOperator{\rank}{rank}
\DeclareMathOperator{\Sel}{Sel}

\newcommand{\lalg}{\lambda^{\alg}}
\newcommand{\malg}{\mu^{\alg}}
\newcommand{\lan}{\lambda^{\an}}
\newcommand{\man}{\mu^{\an}}

\begin{document}

\title{Kida's formula and congruences}
\author{Robert Pollack and Tom Weston}
\address[Robert Pollack]{Department of Mathematics, Boston University, Boston, MA}
\address[Tom Weston]{Dept.\ of Mathematics, University of
Massachusetts, Amherst, MA}

\email[Robert Pollack]{rpollack@math.bu.edu}
\email[Tom Weston]{weston@math.umass.edu}

\thanks{Supported by NSF grants DMS-0439264 and DMS-0440708}

\maketitle

\section{Introduction}

Let $f$ be a modular eigenform of weight at least two and
let $F$ be a finite abelian extension of $\Q$.  Fix an odd prime $p$
at which $f$ is ordinary in the sense that the $p^{\th}$ Fourier
coefficient of $f$ is not divisible by $p$.
In Iwasawa theory, one
associates two objects to $f$ over the cyclotomic $\Zp$-extension
$\Fi$ of $F$: a Selmer group $\Sel(\Fi,A_{f})$ (where $A_{f}$ denotes
the divisible version of the two-dimensional Galois representation
attached to $f$) and a $p$-adic $L$-function $L_{p}(\Fi,f)$.
In this paper we
prove a formula, generalizing work of Kida and
Hachimori--Matsuno, relating the Iwasawa invariants of these objects over $F$
with their Iwasawa invariants over $p$-extensions of $F$.

For Selmer groups our results are significantly more
general.  Let $T$ be a lattice in a 
nearly ordinary $p$-adic Galois representation $V$; set $A=V/T$.
When $\Sel(\Fi,A)$ is a cotorsion Iwasawa module, 
its Iwasawa $\mu$-invariant $\malg(\Fi,A)$ is said to vanish if $\Sel(\Fi,A)$ is cofinitely generated and
its $\lambda$-invariant $\lalg(\Fi,A)$
is simply its $p$-adic corank.
We prove the following result relating these invariants in
a $p$-extension.

\begin{mainthm} \label{mt1}
Let $\Fop/\Fo$ be a finite Galois $p$-extension that is unramified at all places dividing $p$.
Assume that $T$ satisfies the technical assumptions (1)--(5) of Section 2.  
If $\Sel(\Fi,A)$ is $\Lambda$-cotorsion with $\malg(\Fi,A)=0$,
then $\Sel(\Fip,A)$ is $\Lambda$-cotorsion with $\malg(\Fip,A)=0$.
Moreover, in this case
$$\lalg(\Fip,A) = [\Fip:\Fi] \cdot \lalg(\Fi,A)+ 
\sum_{w'} m(\Fipwp/\Fiw,V)$$
where the sum extends over places $w'$ of $\Fip$ which are ramified
in $\Fip/\Fi$.

If $V$ is associated to a cuspform $f$ and  $\Fop$ is an abelian extension of $\Q$, then the same results hold for the analytic Iwasawa invariants of $f$.
\end{mainthm}

Here $m(\Fipwp/\Fiw,V)$ is a certain difference of local multiplicities
defined in Section~\ref{sec:lp}.  In the case of Galois representations
associated to Hilbert modular forms, these local factors can be made quite explicit; see
Section~\ref{sec:hmf} for details.

It follows from Theorem~\ref{mt1} and work of Kato
that if the $p$-adic 
main conjecture holds for a modular form $f$ over $\Q$, then
it holds for $f$ over all abelian $p$-extensions of $\Q$;
see Section~\ref{sec:mc} for details.

These Riemann-Hurwitz type formulas were first discovered by Kida \cite{Kida} in the context of $\lambda$-invariants of CM fields.  More precisely, when $F'/F$ is a $p$-extension of CM fields and  $\mu^-(\Fi/F)=0$, Kida gave a precise formula for $\lambda^-(\Fip/F')$ in terms of $\lambda^-(\Fi/F)$ and local data involving the primes that ramify in $F'/F$.  
(See also \cite{Iwasawa} for a representation theoretic interpretation of
Kida's result.)  This formula was generalized to Selmer groups of
elliptic curves at ordinary primes by Wingberg \cite{Wingberg} in the CM
case and Hachimori--Matsuno \cite{HM} in the general case.
The analytic analogue was first established for ideal class groups by
Sinnott \cite{Sinnott} and for elliptic curves by Matsuno \cite{Matsuno}.

Our proof is most closely related to the arguments in \cite{Sinnott} and \cite{Matsuno} where congruences implicitly played a large role in their study of analytic $\lambda$-invariants.  In this paper, we make the role of congruences more explicit and apply these methods to study both algebraic and analytic $\lambda$-invariants.

As is usual, we first reduce to the
case where $F'/F$ is abelian.  (Some care is required to
show that our local factors are well behaved in towers of fields; this is discussed in Section~\ref{sec:lp}.)
In this case, the $\lambda$-invariant
of $V$ over $F'$ can be expressed as the sum of the $\lambda$-invariants of twists of $V$ by characters of $\Gal(F'/F)$.  The key observation (already visible in both \cite{Sinnott} and \cite{Matsuno}) is that since $\Gal(F'/F)$ is a $p$-group, all of its characters are trivial modulo a prime over $p$ and, thus, the twisted Galois representations  are all congruent to $V$ modulo a prime over $p$.  The algebraic case of
Theorem~\ref{mt1} then follows from the results of
\cite{Weston} which gives a precise local formula for the difference between $\lambda$-invariants of congruent Galois representations.
The analytic case is handled similarly using the results of \cite{EPW}.

The basic principle behind this argument is that a formula relating the Iwasawa invariants of congruent Galois representations should imply of a transition formula for these invariants in $p$-extensions.  As an example of this, in Section \ref{sec:ss}, we use results of \cite{GIP} to prove a Kida formula for the Iwasawa invariants (in the sense of \cite{Perrin-Riou,Kobayashi,Pollack})
of weight 2 modular forms at supersingular primes.

\section{Algebraic invariants}

\subsection{Local preliminaries} \label{sec:lp}

We begin by studying the local terms that appear in our results.
Fix distinct primes $\ell$ and $p$ and let
$L$ denote a finite extension of the cyclotomic $\Zp$-extension of
$\Ql$.  Fix a field $K$ of characteristic zero and a
finite-dimensional $K$-vector space $V$ endowed with
a continuous $K$-linear action of the absolute Galois group $G_L$ of $L$.
Set
$$m_{L}(V) := \dim_{K} \left( V_{I_{L}} \right)^{G_{L}},$$
the multiplicity of the trivial representation in the 
$I_{L}$-coinvariants of $V$.
Note that this multiplicity is invariant under extension of scalars, so that
we can enlarge $K$ as necessary.

Let $L'$  be a finite Galois $p$-extension of $L$.
Note that $L'$ must be cyclic and totally ramified since $L$
contains the $\Zp$-extension of $\Ql$.  
Let
$G$ denote the Galois group of $L'/L$.
Assuming that $K$ contains all $[L':L]$-power roots of unity,
for a character $\chi : G \to K^{\times}$ of $G$, we set
$V_{\chi} = V \otimes_{K} K(\chi)$ with $K(\chi)$ a one-dimensional
$K$-vector space on which $G$ acts via $\chi$.  We define
$$m(L'/L,V) := \sum_{\chi \in G^{\vee}} m_{L}(V) - m_{L}(V_{\chi})$$
where $G^{\vee}$ denotes the $K$-dual of $G$.

The next result shows how these invariants behave in towers of fields.

\begin{lemma} \label{lemma:delta}
Let $L''$ be a finite Galois $p$-extension of $L$ and let $L'$ be a
Galois extension of $L$ contained in $L'$.
Assume that $K$ contains all $[L'':L]$-power roots of unity.
Then
$$m(L''/L,V) = [L'':L'] \cdot m(L'/L,V) +
m(L''/L',V).$$
\end{lemma}
\begin{proof}
Set $G = \Gal(L''/L)$ and $H = \Gal(L''/L')$.
Consider the Galois group $G_{L}/I_{L''}$ 
over $L$ of the maximal unramified extension
of $L''$.  It sits in an exact sequence
\begin{equation} \label{eq:tr}
0 \to G_{L''}/I_{L''} \to G_{L}/I_{L''} \to G \to 0
\end{equation}
which is in fact split since the maximal unramified extensions of both
$L$ and $L''$ are obtained by adjoining all prime-to-$p$ roots of unity.

Fix a character $\chi \in G^{\vee}$.
We compute
\begin{align*}
m_{L}(V_{\chi}) &= \dim_{K} \bigl( (V_{\chi})_{I_{L}}\bigr)^{G_{L}} \\
&= \dim_{K} \Bigl( \bigl( ( (V_{\chi})_{I_{L''}} )_{G} \bigr)^{G_{L''}}
\Bigr)^{G} \\
&=
\dim_{K} \Bigl( \bigl( ( (V_{\chi})_{I_{L''}} )^{G_{L''}} \bigr)_{G}
\Bigr)^{G} \quad\!\! \text{since (\ref{eq:tr}) is split} \\
&= 
\dim_{K} \Bigl( \bigl( (V_{\chi})_{I_{L''}} \bigr)^{G_{L''}} \Bigr)^{G}
\qquad \text{since $G$ is finite cyclic} \\
&=
\dim_{K} \bigl( (V_{I_{L''}} )^{G_{L''}} \otimes \chi \bigr)^{G}
\qquad \text{since $\chi$ is trivial on $G_{L''}$.}
\end{align*}
The lemma thus follows from the following purely group-theoretical
statement applied with $W = (V_{I_{L''}} )^{G_{L''}}$:
for a finite dimensional representation $W$ of
a finite abelian group $G$ over a
field of characteristic zero containing $\mu_{\# G}$, we have
\begin{multline*}
\sum_{\chi \in G^{\vee}} \bigl( \left< W,1 \right>_{G} - \left< W,\chi \right>_{G} \bigr)
= \\
\# H \cdot \sum_{\chi \in (G/H)^{\vee}} \bigl( \left< W,1 \right>_{G} -
\left< W,\chi \right>_{G} \bigr) + \sum_{\chi \in H^{\vee}} \bigl(
\left< W,1 \right>_{H} - \left< W,\chi \right>_{H} \bigr)
\end{multline*}
for any subgroup $H$ of $G$; here $\left< W,\chi \right>_{G}$ 
(resp.\ $\left< W,\chi \right>_{H}$) is the multiplicity of the character
$\chi$ in $W$ regarded as a representation of $G$ (resp.\ $H$).
To prove this, we compute
\begin{align*}
\sum_{\chi \in G^{\vee}} \bigl( &\left< W,1 \right>_{G} - \left< W,\chi \right>_{G} \bigr)\\
 =~~~&\#G \cdot \left< W,1 \right>_{G} - \left< W,\Ind_{1}^{G}1 \right>_{G} \\
=~~~& \#G \cdot \left< W,1 \right>_{G} - \#H \cdot \left<W,\Ind_{H}^{G}1 
\right>_{G}  +\#H \cdot \left<W,\Ind_{H}^{G}1 \right>_{G} -
\left< W,\Ind_{1}^{G} 1\right>_{G} \\
=~~~& \#H \cdot \!\!\!\! \sum_{\chi \in (G/H)^{\vee}} \!\! \bigl(\left<W,1\right>_{G} -
\left<W,\chi \right>_{G} \bigr) +
\sum_{\chi \in H^{\vee}} \left( \left<W,\Ind_{H}^{G}1
\right>_{G} -
\left<W,\Ind_{H}^{G}\chi\right>_{G} \right) \\
=~~~& \#H\cdot \!\!\!\! \sum_{\chi \in (G/H)^{\vee}} \!\! \bigl(\left<W,1\right>_{G} -
\left<W,\chi \right>_{G} \bigr)
  + \sum_{\chi \in H^{\vee}} \bigl(\left<W,1\right>_{H}
- \left<W,\chi\right>_{H}\bigr)
\end{align*}
by Frobenius reciprocity.
\end{proof}

\subsection{Global preliminaries}\label{sec:gp}

Fix a number field $\Fo$; for simplicity we assume that $\Fo$ is
either totally real or totally imaginary.  Fix also an odd prime $p$
and a finite extension $K$ of $\Qp$; we write $\O$ for the ring of integers
of $K$, $\pi$ for a fixed choice of uniformizer of $\O$, and
$k=\O/\pi$ for the residue field of $\O$.

Let $T$ be a nearly ordinary Galois
representation over $\Fo$ with coefficients in $\O$; that is,
$T$ is a free $\O$-module of some rank $n$ endowed with an $\O$-linear
action of the absolute Galois group $\GFo$, together with a
choice for each place $v$ of $\Fo$ dividing $p$ of a complete flag
$$0 = T_{v}^{0} \subset T_{v}^{1} \subset \cdots \subset T_{v}^{n}=T$$
stable under the action of the decomposition group $\GFov \subseteq \GFo$ of $v$.
We make the following assumptions on $T$:
\begin{enumerate}
\item For each place $v$ dividing $p$
we have
$$\left(T_{v}^{i}/T_{v}^{i-1}\right) \otimes k \not\cong
\left(T_{v}^{j}/T_{v}^{j-1}\right) \otimes k$$ 
as $k[\GFov]$-modules for all $i \neq j$;
\item If $\Fo$ is totally real, then
$\rank T^{c_{v}=1}$ is independent of the archimedean place $v$ (here
$c_{v}$ is a complex conjugation at $v$);
\item If $\Fo$ is totally imaginary, then $n$ is even.
\end{enumerate}

\begin{remark}
The conditions above are significantly more restrictive then are actually
required to apply the results of \cite{Weston}.  As our main interest is
in abelian (and thus necessarily Galois) extensions of $\Q$, we have
chosen to include the assumptions (2) and (3) to simply the exposition.
The assumption (1) is also stronger then necessary: all that is actually
needed is that the centralizer of $T \otimes k$ consists entirely of
scalars and that $\mathfrak{gl}_{n}/\mathfrak{b}_{v}$ has trivial adjoint
$G_{v}$-invariants for all places $v$ dividing $p$; here
$\mathfrak{gl}_{n}$ denotes the $p$-adic Lie algebra of $\GL_{n}$ and
$\mathfrak{b}_{v}$ denotes the $p$-adic Lie algebra of the Borel subgroup
associated to the complete flag at $v$.  In particular, when
$T$ has rank $2$, we may still allow the case that $T \otimes k$ has
the form
$$\left( \begin{array}{cc} \chi & * \\ 0 & \chi \end{array} \right)$$
so long as $*$ is non-trivial.  (Equivalently, if $T$ is associated to
a modular form $f$, the required assumption is that $f$ is
{\it $p$-distinguished}.)
\end{remark}

Set $A = T \otimes_{\O} K/\O$; it is a cofree $\O$-module of corank $n$
with an $\O$-linear action of $\GFo$.  Let $c$ equal the rank of 
$T_{v}^{c_{v}=1}$
(resp.\ $n/2$) if $\Fo$ is totally real (resp.\ totally imaginary) and set
$$A_{v}^{\crit} := \im \bigl( T_{v}^{c} \inj T \surj A \bigr).$$
We define the Selmer group of $A$ over the cyclotomic $\Zp$-extension $\Fi$ of $\Fo$ by
$$\Sel(\Fi,A) = \ker \left( H^{1}(\Fi,A) \to
\left( \underset{w \nmid p}{\oplus} H^{1}(\Fiw,A) \right) \times 
\left( \underset{w \mid p}{\oplus} H^{1}(\Fiw,A/A_{v}^{\crit}) 
\right) \right).$$
The Selmer group $\Sel(\Fi,A)$ is naturally a module for the
Iwasawa algebra $\Lambda_{\O} := \O[[\Gal(\Fi/\Fo)]]$.  If
$\Sel(\Fi,A)$ is $\Lambda_{\O}$-cotorsion (that is, if the dual of
$\Sel(\Fi,A)$ is a torsion $\Lambda_{\O}$-module), then we write
$\malg(\Fi,A)$ and $\lalg(\Fi,A)$ for its Iwasawa invariants; in particular,
$\malg(\Fi,A)=0$ if and only if $\Sel(\Fi,A)$ is a cofinitely generated
$\O$-module, while $\lalg(\Fi,A)$ is the $\O$-corank of $\Sel(\Fi,A)$.

\begin{remark}
In the case that $T$ is in fact an {\it ordinary} Galois representation
(meaning that the action of inertia on each $T_{v}^{i}/T_{v}^{i-1}$
is by an integer power $e_{i}$ (independent of $v$)
of the cyclotomic character such that
$e_{1} > e_{2} > \ldots > e_{n}$), then our Selmer group $\Sel(\Fi,A)$
is simply the Selmer group in the sense of Greenberg of a twist of $A$;
see \cite[Section 1.3]{Weston} for details.
\end{remark}

\subsection{Extensions}
Let $\Fop$ be a finite Galois extension of $\Fo$ with degree equal to a power of $p$.
We write $\Fip$ for the cyclotomic
$\Zp$-extension of $\Fop$ and set $G = \Gal(\Fip/\Fi)$.
Note that $T$ satisfies hypotheses (1)--(3) over $\Fop$ as well, 
so that we
may define $\Sel(\Fip,A)$ analogously to $\Sel(\Fi,A)$.  (For (1) this
follows from the fact that $G_{v}$ acts on $(T^{i}_{v}/T^{i-1}_{v}) \otimes k$ by
a character of prime-to-$p$ order; for (2) and (3) it follows from the fact
that $p$ is assumed to be odd.)

\begin{lemma} \label{lemma:twist}
The restriction map
\begin{equation} \label{eq:selchar}
\Sel(\Fi,A) \to \Sel(\Fip,A)^{G}
\end{equation}
has finite kernel and cokernel.
\end{lemma}
\begin{proof}
This is straightforward from the definitions and the fact that
$G$ is finite and $A$ is cofinitely generated;
see \cite[Lemma 3.3]{HM} for details.
\end{proof}

We can use Lemma~\ref{lemma:twist} to relate the $\mu$-invariants of
$A$ over $\Fi$ and $\Fip$.

\begin{corollary} \label{cor:mu}
If $\Sel(\Fi,A)$ is $\Lambda$-cotorsion with $\malg(\Fi,A)=0$, then
$\Sel(\Fip,A)$ is $\Lambda$-cotorsion with $\malg(\Fip,A)=0$.
\end{corollary}
\begin{proof}
This is a straightforward argument using Lemma~\ref{lemma:twist} 
and Nakayama's lemma for compact local rings;
see \cite[Corollary 3.4]{HM} for details.
\end{proof}

Fix a finite extension $K'$ of $K$ containing all $[F':F]$-power roots
of unity.  Consider a character
$\chi : G \to \O'{}^{\times}$
taking values in the ring
of integers $\O'$ of $K'$; note that $\chi$ is
necessarily even since $[F':F]$ is odd.  We set
$$\Ac = A \otimes_{\O} \O'(\chi)$$
where $\O'(\chi)$ is a free $\O'$-module of rank one with $\GF$-action
given by $\chi$.
If we give $\Ac$ the induced complete flags at
places dividing $p$, then $\Ac$ satisfies hypotheses (1)--(3) and we have
$$A_{\chi,v}^{\crit} = A_{v}^{\crit} \otimes_{\O} \O'(\chi) \subseteq
\Ac$$
for each place $v$ dividing $p$.
We write $\Sel(\Fi,\Ac)$ for the corresponding Selmer group,
regarded as a $\Lambda_{\O'}$-module;
in particular, by $\lalg(\Fi,\Ac)$ we mean the
$\O'$-corank of $\Sel(\Fi,\Ac)$, rather than the $\O$-corank.
We write $G^{\vee}$ for the set of all characters
$\chi : G \to \O'{}^{\times}$.

Note that as $\O'[[G_{F'}]]$-modules we have
$$A \otimes_{\O} \O' \cong \Ac$$
from which it follows easily that
\begin{equation} \label{eq:sel}
\bigl(\Sel(\Fip,A) \otimes_{\O} \O'(\chi)\bigr)^{G} =
\Sel(\Fip,\Ac)^{G}
\end{equation}
Moreover, in the case that $G$ is {\it abelian},
\begin{equation}\label{eqn:decomp}\Sel(\Fip,A) \otimes_{\O} \O' \cong
\oplus_{\chi \in G^{\vee}}
\bigl( \Sel(\Fip,A) \otimes_{\O} \O'(\chi) \bigr)^{G}.
\end{equation}
Applying Lemma~\ref{lemma:twist} to each twist $\Ac$, we
obtain the following decomposition of $\Sel(\Fip,A)$.

\begin{corollary}
Assume that $G$ is an abelian group.  Then the map
$$\underset{\chi \in G^{\vee}}{\oplus}
\Sel(\Fi,\Ac) \to \Sel(\Fip,A) \otimes_{\O} \O'$$
obtained from the maps (\ref{eq:selchar}),
(\ref{eq:sel}) and (\ref{eqn:decomp})
has finite kernel and cokernel.
\end{corollary}

As an immediate corollary, we have the following.

\begin{corollary} \label{cor:lambda}
If  $\Sel(\Fi,A)$ is $\Lambda$-cotorsion with $\malg(\Fi,A)=0$, then each group
$\Sel(\Fi,\Ac)$ is $\Lambda_{\O'}$-cotorsion with $\malg(\Fi,\Ac)=0$.  Moreover, if $G$ is abelian, then 
$$\lalg(\Fip,A) = \sum_{\chi \in G^{\vee}} \lalg(\Fi,\Ac).$$
\end{corollary}

\subsection{Algebraic transition formula} \label{sec:statementalg}

We continue with the notation of the previous section.  We write
$R(\Fip/\Fi)$ for the set of prime-to-$p$ places of $\Fip$ which are
ramified in $\Fip/\Fi$.  For a place $w' \in R(\Fip/\Fi)$, we write
$w$ for its restriction to $\Fi$.

\begin{theorem} \label{thm:alg}
Let $\Fop/\Fo$ be a finite Galois $p$-extension with Galois group $G$ which is
unramified at all places dividing $p$.  Let $T$ be a nearly ordinary
Galois representation over $F$ with coefficients in $\O$ satisfying
(1)--(3).  Set $A = T \otimes K/\O$ and assume that:
\begin{enumerate}
\setcounter{enumi}{3}
\item $H^{0}(\Fo,A[\pi]) = H^{0}\bigl(\Fo,\Hom(A[\pi],\mu_{p})\bigr)=0$;
\item $H^{0}(\Iv,A/A_{v}^{\crit})$ is $\O$-divisible for all
$v$ dividing $p$.
\end{enumerate}
If $\Sel(\Fi,A)$ is $\Lambda$-cotorsion with $\malg(\Fi,A)=0$,
then $\Sel(\Fip,A)$ is $\Lambda$-cotorsion with $\malg(\Fip,A)=0$.
Moreover, in this case,
$$\lalg(\Fip,A) = [\Fip:\Fi] \cdot \lalg(\Fi,A)+ 
\sum_{w' \in R(\Fip/\Fi)} \!\!\!\! m(\Fipwp/\Fiw,V)$$
with $V = T \otimes K$ and $m(\Fipwp/\Fiw,V)$ as in Section~\ref{sec:lp}.
\end{theorem}

Note that $m(\Fipwp/\Fiw,V)$ in fact depends only on $w$ and not on $w'$.
The hypotheses (4) and (5) are needed to apply the results of
\cite{Weston}; they will not otherwise appear in the proof below.
We note that the assumption that $\Fop/\Fo$ is unramified at $p$ is
primarily needed to assure that the condition (5) holds for
twists of $A$ as well.

Since $p$-groups are solvable and the only simple $p$-group  is cyclic,
the next lemma shows that it suffices to consider the case of
$\Z/p\Z$-extensions.

\begin{lemma} \label{lemma:reduce}
Let $\Fopp/\Fo$ be a Galois $p$-extension of number fields and let $\Fop$ be an intermediate extension which is Galois over $\Fo$.  Let $T$ be as above.
If Theorem \ref{thm:alg} holds for $T$ with respect to any two of the three field extensions $\Fopp/\Fop$,  $\Fop/\Fo$ and $\Fopp/\Fo$, then it holds  for $T$ with respect to the third extension.
\end{lemma}

\begin{proof}
This is clear from Corollary~\ref{cor:mu} except for the
$\lambda$-invariant formula.  Substituting the formula for
$\lambda(\Fip,A)$ in terms of $\lambda(\Fi,A)$ into the formula for
$\lambda(\Fipp,A)$ in terms of $\lambda(\Fip,A)$, one finds that it
suffices to show that
\begin{multline*}
\sum_{w'' \in R(\Fipp/\Fi)} m(\Fippwpp/\Fiw,V) = \\
[\Fipp:\Fip] \cdot \sum_{w' \in R(\Fip/\Fi)} m(\Fipwp/\Fiw,V) \\ + \sum_{w'' \in R(\Fipp/\Fip)} m(\Fippwpp/\Fipwp,V).
\end{multline*}
This formula follows upon summing the formula of Lemma~\ref{lemma:delta}
over all $w'' \in R(\Fipp/\Fi)$ and using the two facts:
\begin{itemize}
\item $[\Fipp:\Fip]/[\Fippwpp:\Fipwp]$ equals the number of places of
$\Fipp$ lying over $w'$ (since 
the residue field of $\Fiw$ has no $p$-extensions);
\item $m(\Fippwpp/\Fipwp,V) = 0$ for any $w'' \in R(\Fipp/\Fi)-R(\Fipp/\Fip)$.
\end{itemize}
\end{proof}

\begin{proof}[Proof of Theorem \ref{thm:alg}]
By Lemma~\ref{lemma:reduce} and the preceding remark, we may assume
that $\Fip/\Fi$ is a cyclic extension of degree $p$.
The fact that $\Sel(\Fip,A)$ is cotorsion with trivial $\mu$-invariant is
simply Corollary~\ref{cor:mu}.
Furthermore, by Corollary~\ref{cor:lambda}, we have
$$\lalg(\Fip,A) = \sum_{\chi \in G^{\vee}}
\lalg(\Fi,\Ac).$$

For $\chi \in G^{\vee}$, note that $\chi$ is trivial modulo a uniformizer
$\pi'$ of $\O'$
as it takes values in $\mu_{p}$.  In particular, the residual
representations $\Ac[\pi']$ and $A[\pi]$ are isomorphic.  
Under the hypotheses (1)--(5), the result
\cite[Theorem 1]{Weston} gives a precise formula for the relation between
$\lambda$-invariants of congruent Galois representations.  In the
present case it takes the form:
$$\lalg(\Fi,\Ac) = \lalg(\Fi,A) + \sum_{w' \nmid p} \bigl( m_{\Fiw}(V \otimes \omega^{-1}) - 
m_{\Fiw}(V_{\chi} \otimes \omega^{-1}) \bigr)$$
where the sum is over all prime-to-$p$ places $w'$ of $\Fip$,
$w$ denotes the place of $\Fi$ lying under $w'$ and $\omega$ is the mod $p$ cyclotomic character.
The only non-zero terms in this sum are those for which $w'$ is ramified
in $\Fip/\Fi$.  For any such $w'$, we have $\mu_{p} \subseteq \Fiw$ by local
class field theory so that
$\omega$ is in fact trivial at $w$; thus
$$\lalg(\Fi,\Ac) = \lalg(\Fi,A) + \sum_{w' \in R(\Fip/\Fi)}\bigl( m_{\Fiw}(V) - m_{\Fiw}(V_{\chi}) \bigr).$$
Summing over all $\chi \in G^\vee$ then yields
$$
\lalg(\Fip,A) = [\Fip:\Fi] \cdot \lalg(\Fi,A) + \sum_{w' \in R(\Fip/\Fi)} \!\!\!\!
m(\Fipwp/\Fiw,V)$$
which completes the proof. 
\end{proof}

\section{Analytic invariants}

\subsection{Definitions} \label{sec:def}
Let $f = \sum a_{n}q^{n}$ be a modular eigenform of weight $k \geq 2$, level $N$
and character $\varepsilon$.
Let $K$ denote the finite extension of $\Qp$ generated by
the Fourier coefficients of $f$ (under some fixed embedding
$\Qbar \inj \Qpbar$), let $\O$ denote the ring of integers of $K$ and
let $k$ denote the residue field of $\O$.
Let $V_{f}$ denote a two-dimensional $K$-vector space with Galois action
associated to $f$ in the usual way; thus the characteristic polynomial of
a Frobenius element at a prime $\ell \nmid Np$ is
$$x^{2} - a_{\ell}x + \ell^{k-1}\varepsilon(\ell).$$
Fix a Galois stable $\O$-lattice $T_{f}$ in $V_{f}$.  We assume that
$T_{f} \otimes k$ is an irreducible Galois representation; in this case
$T_{f}$ is uniquely determined up to scaling.  Set
$A_{f} = T_{f} \otimes K/\O$.

Assuming that $f$ is $p$-ordinary (in the sense that $a_{p}$ is relatively
prime to $p$) and fixing a canonical period for $f$, one can associate to $f$ a
$p$-adic $L$-function $L_{p}(\Qinf/\Q,f)$ which lies in  $\Lambda_\O$.
This is
well-defined up to a $p$-adic unit (depending upon the choice  of a canonical period) and thus has well-defined Iwasawa invariants.

Let $\Fo/\Q$ be a finite abelian extension
and let $\Fi$ denote the cyclotomic $\Zp$-extension of $\Fo$.  For a character $\chi$ of $\Gal(\Fo/\Q)$, we denote by $f_{\chi}$ the modular eigenform $\sum a_{n}\chi(n)q^{n}$
obtained from $f$ by twisting by $\chi$ (viewed as a Dirichlet character).  If $f$ is $p$-ordinary and $\Fo/\Q$ is unramified at $p$, then $f_\chi$ is again $p$-ordinary and we define
$$
L_p(\Fi/\Fo,f) = \prod_{\chi \in \Gal(\Fo/\Q)^\vee} L_p(\Qinf/\Q,f_\chi).
$$
If $\Fo/\Q$ is ramified at $p$, it is still possible to
 define $L_p(\Fi/\Fo,f)$; see \cite[pg.\ 5]{Matsuno}, for example.


If $F_1$ and $F_2$ are two distinct number fields whose cyclotomic $\Zp$-extensions agree,
the corresponding $p$-adic $L$-functions of $f$ over $F_1$ and $F_2$ need not agree.  However, it is easy to check that the Iwasawa invariants of these two power series are equal.  We thus denote the Iwasawa invariants of $L_p(\Fi/F,f)$ simply by $\man(\Fi,f)$ and $\lan(\Fi,f)$.

\subsection{Analytic transition formula}
\label{sec:statementan}

Let $F/\Q$ be a finite abelian $p$-extension of $\Q$  and
let $\Fop$ be a finite $p$-extension of $\Fo$ such that $\Fop/\Q$ is abelian.  As always, let $\Fi$ and $\Fip$ denote the cyclotomic $\Zp$-extensions of $\Fo$ and $\Fop$.
As before, we write $R(\Fip/\Fi)$ for the set of prime-to-$p$ places of $\Fip$
which are ramified in $\Fip/\Fi$.

\begin{theorem} \label{thm:an}
Let $f$ be a $p$-ordinary modular form such that $T_f \otimes k$ is irreducible and $p$-distinguished.
If $\man(\Fi,f)=0$, then $\man(\Fip,f)=0$.  Moreover, if this is the case, then 
$$
\label{eqn:lan}
\lan(\Fip,f) = [\Fip:\Fi] \cdot \lan(\Fi,f) + 
\sum_{w' \in R(\Fip/\Fi)} \!\!\!\!m(\Fipwp/\Fiw,V_{f}).
$$
\end{theorem}

\begin{proof}
By Lemma~\ref{lemma:reduce}, we may assume $[F:\Q]$ is  prime-to-$p$.  Indeed,
let $F_0$ be the maximal  subfield of $F$ of prime-to-$p$ degree over $\Q$.
By Lemma~\ref{lemma:reduce}, knowledge of the theorem for the two
extensions $F'/F_{0}$ and $F/F_{0}$ would then imply it for $F'/F$ as well.

We may further assume that $F$ and $F'$ are unramified at $p$.  Indeed, if $F^{\text{ur}}$ (resp.\ $F'{}^{\text{ur}}$) denotes the maximal subfield of $\Fi$ (resp.\ $\Fip$) unramified at $p$, then $F^{\text{ur}} \subseteq F'{}^{\text{ur}}$ and
the cyclotomic $\Zp$-extension of $F^{\text{ur}}$ (resp.\ $F'{}^{\text{ur}}$) 
is $\Fi$ (resp.\ $\Fip$).  
Thus, by the comments at the end of Section \ref{sec:def}, we may replace $F$ by $F^{\text{ur}}$ and $F'$ by $F'^{\text{ur}}$ without altering the formula we are studying.

After making these reductions, we let $M$ denote the (unique) $p$-extension of $\Q$ inside of $\Fop$ such that $M\Fo = \Fop$.  Set $G = \Gal(F/\Q)$ and $H = \Gal(M/\Q)$, so that $\Gal(F'/\Q) \cong G \times H$.  Then since $F$ and $F'$ are
unramified at $p$ by definition, we have
\begin{equation} \label{eq:man1}
\man(\Fi,f) = \sum_{\psi \in \Gal(\Fo/\Q)^\vee} \man(\Qinf,f_\psi)
\end{equation} 
and
\begin{equation} \label{eq:man2}
\man(\Fip,f) = \sum_{\psi \in \Gal(\Fop/\Q)^\vee} \man(\Qinf,f_\psi) = \sum_{\psi \in G^{\vee}} \sum_{\chi \in H^{\vee}} 
\man(\Qinf,f_{\psi\chi}).
\end{equation}
Since we are assuming that $\man(\Fi,f) = 0$ and since these $\mu$-invariants are non-negative, from (\ref{eq:man1}) it follows that $\man(\Qinf,f_\psi) = 0$ for each $\psi \in \Gal(\Fo/\Q)^\vee$.

Fix $\psi \in G^{\vee}$.
For any $\chi \in H^\vee$, $\psi\chi$ is congruent to $\psi$ 
modulo any prime over $p$ and  thus $f_\chi$ and $f_{\psi\chi}$ are congruent modulo any prime over $p$.  Then, since $\man(\Qinf,f_\psi)=0$, by \cite[Theorem 1]{EPW} it follows that $\man(\Qinf,f_{\psi\chi})=0$ for each $\chi \in H^\vee$.  Therefore, by (\ref{eq:man2}) we have that $\man(\Fip,f)=0$ proving the first part of the theorem.

For $\lambda$-invariants, we again have
$$
\lan(\Fi,f) = \sum_{\psi \in \Gal(\Fo/\Q)^\vee} \lan(\Qinf,f_\psi).
$$
and
\begin{equation} \label{eq:lan}
\lan(\Fip,f) = \sum_{\psi \in G^{\vee}} \sum_{\chi \in H^{\vee}} 
\lan(\Qinf,f_{\psi\chi}).
\end{equation}
By \cite[Theorem 2]{EPW} the congruence between
$f_\chi$ and $f_{\psi\chi}$ implies that
\begin{multline*}
\lan(\Qinf,f_{\psi\chi}) - \lan(\Qinf,f_{\psi}) =\\
\sum_{v' \in R(M_{\infty}/\Qinf)} \!\!\!\!\!\left(
m_{\Qinfv}(V_{f_{\psi\chi}} \otimes \omega^{-1}) 
- m_{\Qinfv}(V_{f_{\psi}} \otimes \omega^{-1}) \right)
\end{multline*}
where $v$ denotes the place of $\Qinf$ lying under the place $v'$ of
$M_{\infty}$.
Note that in \cite{EPW} the sum extends over all prime-to-$p$ places;
however, the terms are trivial unless $\chi$ is ramified at $v$.  Also note that
the mod $p$ cyclotomic characters that appear are actually trivial since if $\Qinfv$ has a ramified Galois $p$-extensions for $v \nmid p$, then $\mu_p \subseteq \Qinfv$.

Combining this with (\ref{eq:lan}) and the definition of
$m(M_{\infty,v'}/\Qinfv,V_{f_{\psi}})$, we conclude that
\begin{align*}
\lan(\Fip,f) &= \sum_{\psi \in G^{\vee}} 
\left( [\Fip:\Fi] \cdot \lan(\Qinf,f_{\psi}) + \sum_{v' \in R(M_\infty/\Qinf)} \!\!\!\!\!
m(M_{\infty,v'}/\Qinfv,V_{f_{\psi}}) \right) \\
&= [\Fip:\Fi] \cdot \lan(\Fi,f) + \sum_{v' \in R(M_\infty/\Qinf)} \sum_{\psi \in G^{\vee}} 
m(M_{\infty,v'}/\Qinfv,V_{f_{\psi}}) \\
&= [\Fip:\Fi] \cdot \lan(\Fi,f) + \sum_{v' \in R(M_\infty/\Qinf)} 
g_{v'}(\Fip/M_{\infty}) \cdot \\ &
\qquad \qquad m(M_{\infty,v'}/\Qinfv,\Z[\Gal(\Fiw/\Qinfv)] \otimes V_{f})
\end{align*}
where $g_{v'}(\Fip/M_{\infty})$ denotes the number of places of
$\Fip$ above the place $v'$ of $M_{\infty}$.
By Frobenius reciprocity,
$$m(M_{\infty,v'}/\Qinfv,\Z[\Gal(\Fiw/\Qinfv)] \otimes V_{f}) =
m(\Fipwp/\Fiw,V_{f})$$
where $w'$ is the unique place of $\Fip$ above $v'$ and $w$. It follows that
$$\lambda(\Fip,f) = [\Fip:\Fi] \cdot \lan(\Fi,f) + \sum_{w' \in R(\Fip/\Fi)}  
m(\Fipwp/\Fiw,V_{f})$$
as desired.
\end{proof}

\section{Additional Results}

\subsection{Hilbert modular forms} \label{sec:hmf}

We illustrate our results in the case of the two-dimensional
representation $V_{f}$ associated to a Hilbert modular eigenform $f$ 
over a totally real field
$\Fo$.  Although in principle our analytic results should remain true
in this context, we focus on the less conjectural algebraic picture.
Fix a $\GFo$-stable lattice $T_{f} \subseteq V_{f}$
and let $A_{f} = T_{f} \otimes K/\O$.

Let $\Fop$ be a finite Galois $p$-extension of $\Fo$
unramified at all places dividing $p$; for simplicity we assume also
that $\Fop$ is linearly disjoint from $\Fi$.
Let $v$ be a place of $\Fo$ not dividing $p$ and fix a place $v'$ of $\Fop$ lying over $v$.
For a character $\varphi$ of $\GFov$, we define
$$h(\varphi) = \begin{cases}
-1 & \varphi \text{~ramified,~} \varphi|_{\GFopvp} 
\text{~unramified, and~}
\varphi \equiv 1 \bmod{\pi}; \\
0 & \varphi \not\equiv 1 \bmod{\pi} \text{~or~} \varphi|_{\GFopvp} 
\text{~ramified}; \\
e_{v}(\Fop/\Fo) - 1 & \varphi \text{~unramified and~} \varphi \equiv 1
\bmod{\pi} \end{cases}$$
where $e_{v}(\Fop/\Fo)$ denotes the ramification index of $v$ in $\Fop/\Fo$ and $\GFopvp$ is the decomposition group at $v'$. Set
$$h_{v}(f) = \begin{cases}
h(\varphi_{1}) + h(\varphi_{2}) &
f \text{~principal series with characters~} \varphi_{1},\varphi_{2}
\text{~at~} v; \\
h(\varphi) & f \text{~special with character~} \varphi \text{~at~} v; \\
0 & f \text{~supercuspidal or extraordinary at~} v.
\end{cases}$$
For example, if $f$ is unramified principal series at $v$ with
Frobenius characteristic polynomial
$$x^{2}-a_{v}x + c_{v},$$
then
$$h_{v}(f) = \begin{cases} 2(e_{v}(\Fop/\Fo)-1) &
a_{v} \equiv 2, c_{v} \equiv 1 \bmod{\pi} \\
e_{v}(\Fop/\Fo)-1 & a_{v} \equiv c_{v}+1 \not\equiv 2 \bmod{\pi} \\
0 & \text{otherwise}.
\end{cases}
$$

\begin{theorem} \label{thm:hmf}
Assume that $f$ is ordinary (in the sense that for each place $v$ dividing $p$
the Galois representation
$V_{f}$ has a unique one-dimensional quotient unramified at $v$) and that
$$H^{0}(\Fo,A_{f}[\pi]) = H^{0}\bigl(\Fo,\Hom(A_{f}[\pi],\mu_{p})\bigr)=0
.$$
If $\Sel(\Fi,A_{f})$ is $\Lambda$-cotorsion with $\malg(\Fi,A_{f})=0$,
then also $\Sel(\Fip,A_{f})$ is $\Lambda$-cotorsion with
$\malg(\Fip,A_{f})=0$ and
$$\lalg(\Fip,A) = [\Fip:\Fi] \cdot \lalg(\Fi,A)+
\sum_{v}
g_{v}(\Fip/\Fo) \cdot h_{v}(f);$$
here the sum is over the prime-to-$p$ places of $F$ ramified in $\Fip$ and
$g_{v}(\Fip/\Fo)$ denotes the number of places of $\Fip$ lying over such a $v$.
\end{theorem}
\begin{proof}
Fix a place $v$ of $\Fo$ not dividing $p$ and let $w$ denote
a place of $\Fi$ lying over $v$.
Since there are exactly $g_{v}(\Fi/\Fo)$ such places,
by Theorem~\ref{thm:alg} it suffices to prove that
\begin{multline} \label{eq:hvf}
h_{v}(f) = m(\Fipwp/\Fiw,V_{f}) := \\
\sum_{\chi \in \Gal(\Fipwp/\Fiw)^{\vee}} \left( m_{\Fiw}(V_{f}) -
m_{\Fiw}(V_{f,\chi}) \right).
\end{multline}
This is a straightforward case
analysis.  We will discuss the case that $V_{f}$ is special associated
to a character $\varphi$ at $v$; the other cases are similar.
In the special case, we have
$$V_{f,\chi}|_{I_{\Fiw}} = \begin{cases}
K'(\chi\varphi) & \chi\varphi|_{G_{\Fiw}} \text{~unramified}; \\
0 & \chi\varphi|_{G_{\Fiw}} \text{~ramified}.
\end{cases}$$
Since an unramified character has trivial restriction to $G_{\Fiw}$
if and only if it has trivial reduction modulo $\pi$, it follows that
$$m_{\Fiw}(V_{f,\chi}) = \begin{cases}
1 & \varphi \equiv 1 \bmod{\pi} \text{~and~} \chi\varphi|_{G_{\Fiw}} 
\text{~unramified}; \\
0 & \text{otherwise}.
\end{cases}$$
In particular, the sum in (\ref{eq:hvf}) is zero if
$\varphi \not\equiv 1 \bmod{\pi}$
or if $\varphi$ is ramified when restricted to
$G_{\Fipwp}$ 
(as then $\chi\varphi$ is ramified for all $\chi \in \Gv^{\vee}$).
If $\varphi \equiv 1 \bmod{\pi}$ and
$\varphi$ itself is unramified, then $m_{\Fiw}(V_{f})=1$ while
$m_{\Fiw}(V_{f,\chi}) = 0$ for $\chi \neq 1$, so that the sum
in (\ref{eq:hvf}) is $[\Fipwp:\Fiw]-1=e_{v}(\Fop/\Fo)-1$,
as desired.  Finally,
if $\varphi  \equiv 1 \bmod{\pi}$ and $\varphi$ is ramified but
becomes unramified when restricted to $\GFopvp$, then $m_{\Fiw}(V_{f})=0$,
while $m_{\Fiw}(V_{f,\chi})=1$ for a unique $\chi$, so that the sum is $-1$.
\end{proof}

Suppose finally that $f$ is in fact the Hilbert modular form associated
to an elliptic curve $E$ over $\Fo$.  The only principal series
which occur are unramified and we have $c_{v} \equiv 1 \pmod{\pi}$
(since the determinant of $V_{f}$ is cyclotomic and $\Fi$ has
a $p$-extension (namely, $\Fip$) ramified at $v$),
 so that
$$h_{v}(f) \neq 0 \quad \Leftrightarrow \quad a_{v} =2 \quad 
\Leftrightarrow \quad
E(\Fov) \text{~has a point of order $p$}$$
in which case $h_{v}(f) = 2(e_{v}(\Fop/\Fo)-1)$.  The only characters which
may occur in a special constituent are trivial or unramified quadratic,
and we have $h_{v}(f)=e_{v}(\Fop/\Fo)-1$ or $0$ respectively.  Thus
Theorem~\ref{thm:hmf} recovers \cite[Theorem 3.1]{HM} in this case.

\subsection{The main conjecture} \label{sec:mc}

Let $f$ be a $p$-ordinary elliptic
modular eigenform of weight at least two and arbitrary
level with associated Galois representation $V_{f}$.  
Let $F$ be a finite abelian extension of $\Q$ with cyclotomic
$\Zp$-extension $\Fi$.
Recall that the $p$-adic Iwasawa main conjecture for $f$ over $F$
asserts that the
Selmer group $\Sel(\Fi,A_{f})$ is $\Lambda$-cotorsion and that the
characteristic ideal of its dual
is generated by the $p$-adic $L$-function
$L_{p}(\Fi,f)$.  In fact, when the residual representation of $V_{f}$
is absolutely irreducible, it is known by work of Kato that
$\Sel(\Fi,A_{f})$ is indeed $\Lambda$-cotorsion and that
$L_{p}(\Fi,f)$ is an element of the characteristic ideal of
$\Sel(\Fi,A_{f})$.  In particular, this reduces the verification of
the main conjecture for $f$ over $F$ to the equality of the algebraic
and analytic Iwasawa invariants of $f$ over $F$.
The identical transition formulae in
Theorems~\ref{thm:alg} and~\ref{thm:an} thus yield the following immediate
application to the main conjecture.

\begin{theorem} \label{thm:mc}
Let $F'/F$ be a finite $p$-extension with $F'$ abelian over $\Q$.
If the residual representation of $V_{f}$ is absolutely irreducible and $p$-distinguished,
then the main conjecture holds for $f$ over $F$ with $\mu(\Fi,f)=0$ if and only if it holds for
$f$ over $F'$ with $\mu(\Fip,f)=0$ .
\end{theorem}

We note that in Theorem~\ref{thm:alg}, it was assumed that $F'/F$ was unramified at all places over $p$.  However, in this special case where $F'/\Q$ is abelian, this hypothesis can be removed.  Indeed, one simply argues in an analogous way as at the start of Theorem~\ref{thm:an} by replacing $F'$ (resp.\  $F$) by the maximal sub-extension of $\Fip$ (resp.\ $\Fi$) that is unramified at $p$.

For an example of Theorem \ref{thm:mc}, consider the eigenform 
$$\Delta = q \prod_{n \geq 1} (1-q^{n})^{24}$$ 
of weight $12$ and level $1$.  We take $p=11$.  It is well known that
$\Delta$ is congruent modulo $11$ to the newform associated to the
elliptic curve $X_{0}(11)$.  The $11$-adic main conjecture is known for 
$X_{0}(11)$ over $\Q$; it has trivial $\mu$-invariant and
$\lambda$-invariant equal to $1$ (see, for instance, \cite[Example 5.3.1]{EPW}).  We should be clear here that the
non-triviality of $\lambda$ in this case corresponds to a trivial zero of
the $p$-adic $L$-function; we are using the Greenberg Selmer group which
does account for the trivial zero.)  It follows from
\cite{EPW} that the $11$-adic main conjecture also holds for $\Delta$
over $\Q$, again with trivial $\mu$-invariant and $\lambda$-invariant equal
to $1$.  Theorem~\ref{thm:mc} thus allows us to conclude that the main
conjecture holds for $\Delta$ over any abelian $11$-extension of $\Q$.

For a specific example, consider
 $F=\Q(\zeta_{23})^{+}$; it is a cyclic $11$-extension of $\Q$.
We can easily use
Theorem~\ref{thm:hmf} to compute its $\lambda$-invariant:
using that $\tau(23) = 18643272$ one finds that
$h_{23}(\Delta)=0$, so that $\lambda(\Q(\zeta_{23})^{+},\Delta)=11$.

For a more interesting example, take $F$ to be the unique subfield
of $\Q(\zeta_{1123})$ which is cyclic of order $11$ over $\Q$.
In this case we have
$$\tau(1123) \equiv 2 \pmod{11}$$
so that we have $h_{1123}(\Delta) = 20$.  Thus,
in this case, Theorem~\ref{thm:hmf} shows that
$\lambda(F,\Delta) = 31$.

\subsection{The supersingular case}
\label{sec:ss}

As mentioned in the introduction, the underlying principle of this paper is that the existence of a formula relating the $\lambda$-invariants of congruent Galois representations should imply a Kida-type formula for these invariants.  We illustrate this now in the case of modular forms of weight two that are supersingular at $p$.

Let $f$ be an eigenform of weight 2 and level $N$ with Fourier coefficients in $K$ some finite extension of $\Qp$.  Assume further than $p \nmid N$ and that $a_p(f)$ is not a $p$-adic unit.  In \cite{Perrin-Riou}, Perrin-Riou associates
to $f$ a pair of algebraic and analytic $\mu$-invariants over $\Qinf$ which we denote by $\mu^\star_\pm(\Qinf,f)$.  (Here $\star$ denotes either ``$\alg$" or ``$\an$" for algebraic and analytic respectively.)  Moreover, when $\mu^\star_+(\Qinf,f) = \mu^\star_-(\Qinf,f)$ or when $a_p(f)=0$, she also defines corresponding
$\lambda$-invariants $\lambda^\star_\pm(\Qinf,f)$.  When $a_p(f)=0$ these invariants coincide with the Iwasawa invariants of \cite{Kobayashi} and \cite{Pollack}.  We also note that in \cite{Perrin-Riou} only the case of elliptic curves is treated, but the methods used there generalize to weight two modular forms.

We extend the definition of these invariants to the cyclotomic $\Zp$-extension of an abelian extension $F$ of $\Q$.  As usual, by passing to the
maximal subfield of $\Fi$ unramified at $p$, we may assume that $F$ is
unramified at $p$.
We define
$$
\mu^\star_\pm(\Fi,f) = \sum_{\psi \in \Gal(F/\Q)^{\vee}} \!\!\!\! \mu^\star_\pm(\Qinf,f_\psi) \hspace{.4cm} \text{and} \hspace{.4cm}
\lambda^\star_\pm(\Fi,f) = \sum_{\psi \in \Gal(F/\Q)^{\vee}} \!\!\!\! \lambda^\star_\pm(\Qinf,f_\psi)
$$
for $\star \in \{\alg,\an\}$.  

The following transition formula follows from the congruence results of \cite{GIP}.

\begin{theorem} \label{thm:ssan}
Let $f$ be as above and assume further 
that $f$ is congruent modulo some prime above $p$ to a modular form with
coefficients in $\Zp$.
Consider an extension of number fields $F'/F$ with $F'$ an
abelian $p$-extension of $\Q$.
If $\mu^\star_\pm(\Fi,f)=0$, then $\mu^\star_\pm(\Fip,f)=0$.  Moreover, if this is the case, then 
$$
\lambda^\star_\pm(\Fip,f) = [\Fip:\Fi] \cdot \lambda^\star_\pm(\Fi,f) + 
\sum_{w' \in R(\Fip/\Fi)} \!\!\!\! m(\Fipwp/\Fiw,V_{f}).
$$
In particular, if the main conjecture is true for $f$ over $F$ (with $\mu^\star_\pm(\Fi,f)=0$), then the main conjecture is true for $f$ over $F'$ (with $\mu^\star_\pm(\Fip,f)=0$).
\end{theorem}

\begin{proof}
The proof of this theorem proceeds along the lines of the proof of  Theorem \ref{thm:an} replacing the appeals to the results of \cite{EPW,Weston} to the results of \cite{GIP}.  The main result of \cite{GIP} is a formula relating the $\lambda^\star_\pm$-invariants of congruent supersingular weight two modular forms.  This formula has the same shape as the formulas that appear in \cite{EPW} and \cite{Weston} which allows for the proof to proceed nearly verbatim.  The hypothesis that $f$ be congruent to a modular form with $\Zp$-coefficients
is needed because this hypothesis appears in the results of \cite{GIP}.

One difference to note is that in this proof 
we need to assume that $F$ is a $p$-extension of $\Q$.  The reason for this assumption is that in the course of the proof we need to apply the results of \cite{GIP} to the form $f_\psi$ where $\psi \in \Gal(F/\Q)^\vee$. We thus need to know that $f_\psi$ is congruent to some modular form with coefficients in $\Zp$.  In the case that $\Gal(F/\Q)$ is a $p$-group, $f_\psi$ is congruent to $f$ which by assumption is congruent to such a form.
 \end{proof}

\end{document}